\documentclass[12pt]{amsart}

\newtheorem{theo+}              {Theorem}           [section]
\newtheorem{prop+}  [theo+]     {Proposition}
\newtheorem{coro+}  [theo+]     {Corollary}
\newtheorem{lemm+}  [theo+]     {Lemma}
\newtheorem{exam+}  [theo+]     {Example}
\newtheorem{rema+}  [theo+]     {Remark}
\newtheorem{defi+}  [theo+]     {Definition}

\newenvironment{theorem}{\begin{theo+}}{\end{theo+}}
\newenvironment{proposition}{\begin{prop+}}{\end{prop+}}
\newenvironment{corollary}{\begin{coro+}}{\end{coro+}}
\newenvironment{lemma}{\begin{lemm+}}{\end{lemm+}}
\newenvironment{example}{\begin{exam+}}{\end{exam+}}
\newenvironment{remark}{\begin{rema+}}{\end{rema+}}
\newenvironment{definition}{\begin{defi+}}{\end{defi+}}
\theoremstyle{definition}

\newtheorem*{ack}{{\bf Acknowledgments}}

\renewcommand{\Bbb}{\mathbb}

\def\E{/\kern-1.0em \equiv }
\def\h{harmonic morphism}
\def\qh{quadratic harmonic morphism}
\def\om{orthogonal multiplication}

\begin{document}
\baselineskip  18pt
\title [Quadratic harmonic morphisms and O-systems]{Quadratic harmonic
morphisms and O-systems}
\bigskip
\thanks{$^{*}$The author was supported by the Chinese Education  Commission \\
as ``Visiting scholar program (94A04)''.}
\maketitle

\begin{center}
{\bf Ye-Lin Ou$^{*}$}\\
Department of Mathematics, Guangxi University for Nationalities, Nanning
530006, P.R.China.\\
Current : Department of Pure Mathematics, School of Mathematics, University of
Leeds, Leeds LS2 9JT, U.K. (until Jan.1996). \\
E-mail: pmtyo@amsta.leeds.ac.uk\\
\end{center}
\bigskip

\begin{center}
{\bf Abstract}
\end{center}
\bigskip
We introduce O-systems (Definition \ref{DO}) of orthogonal transformations of
${\Bbb R}^{m}$, and establish $1-1$ correspondences both between equivalence
classes of Clifford systems and that of O-systems, and between O-systems and
orthogonal multiplications of the form $\mu :{\Bbb R}^{n} \times {\Bbb R}^{m}
\longrightarrow {\Bbb R}^{m} $, which allow us to solve the existence problems
both for O-systems and for umbilical quadratic harmonic morphisms (Theorems
\ref{ES} and \ref{EU}) simultaneously. The existence problem for general
quadratic harmonic morphisms is then solved (Theorem \ref{EG}) by the Splitting
Lemma (Lemma \ref{Split}). We also study properties (see, e.g., Theorems
\ref{single} and \ref{TL}) possessed by all quadratic harmonic morphisms for
fixed pairs of domain and range spaces (\S5).
\bigskip

{Keywords: Harmonic maps, Quadratic harmonic morphisms,\\~O-systems, Clifford
systems.}
\newpage

\section{Introduction}

A map  $\varphi:({M}^{m},g)\longrightarrow ({N}^{n},h)$ between two Riemannian
manifolds is called a {\bf harmonic map} if the divergence of its differential
vanishes. Such maps are the extrema  of the energy functional $
\frac{1}{2}\int_{D} {\left|d\varphi\right|}^2$ over compact domain D in M. For
a detailed account on harmonic maps we refer to \cite{EelLem78, EelLem83,
EelLem88} and the references therein. {\bf Harmonic morphisms} are a special
subclass of harmonic maps which preserve solutions of Laplace's equation in the
sense that for any harmonic function $f:U \longrightarrow{\Bbb R}$, defined on
an open subset U of $N$ with ${\varphi}^{-1}(U) $ non-empty, $f\circ\varphi
:{\varphi}^{-1}(U) \longrightarrow {\Bbb R} $ is a harmonic function. In other
words, $\varphi$ pulls back germs of harmonic functions on $N$ to germs of
harmonic functions on $M$. In the theory of stochastic process, harmonic
morphisms $\varphi:(M,g)\longrightarrow (N,h)$ are found to be Brownian pat!
 h preserving mappings meaning that they map Brownian motions on $M$ to
Brownian motions on $N $(\cite{BerCamDav79, Levy}). It is well-known (see
\cite{Fug78, Ish79A}) that a map between Riemannian manifolds is a harmonic
morphism if and only if it is both a harmonic map and a horizontally weakly
conformal map.
For a map  $\varphi:{\Bbb  R}^{m}\longrightarrow{\Bbb R}^{n}$ between Euclidean
spaces, with $\varphi(x)=({\varphi}^{1}(x),...,{\varphi}^{n}(x))$, the
harmonicity and horizontally weakly conformality are equivalent to the
following conditions respectively:
\begin{align}
&\sum_{i=1}^{m} \frac{\partial^{2}{\varphi}^{\alpha}}{\partial{x}_{i}^{2}} = 0
\label{1a}\\
&\sum_{i=1}^{m} \frac{\partial{\varphi}^{\alpha}}{\partial{x}_{i}}
\frac{\partial{\varphi}^{\beta}}{\partial{x}_{i}} = {\lambda}^{2}(x)
{\delta}^{\alpha \beta}\ ,\ \alpha, \beta =1,2,...,n.\label{1b}
\end{align}
where $(x_{1},\ldots,x_{m})$ are the standard coordinates of ${\Bbb  R}^{m}$.\\

In recent years, much work has been done in classifying and constructing
harmonic morphisms from certain model spaces to other model spaces (see e.g.
\cite{Bai83, BaiWoo88, BaiWoo91, BaiWoo92A, BaiWoo95, GudA, Gud92, Gud94,
GudSig93, Woo86A, Woo92, OuWoo95, Ou95A, Ou95C}). Concerning \h s between
Euclidean spaces, Baird \cite{Bai83} has studied \h s $\varphi:{\Bbb
R}^{m}\longrightarrow{\Bbb R}^{n}$ defined by homogeneous polynomials of degree
$p$. He has obtained a necessary condition on the dimensions of the domain and
the range spaces for such \h s to exist; he also gives a possible way to
construct such \h s from a single polynomial. For \qh s, i.e., \h s defined by
homogeneous polynomials of degree $2$, he proves ( Theorem 7.2.7 in
\cite{Bai83}) that an orthogonal multiplication $ \mu : {\Bbb R}^{p} \times
{\Bbb R}^{q} \longrightarrow {\Bbb R}^{n}$ is a harmonic morphism if and only
if the dimensions $p = q = n$ and $n = 1,2,4,$ or $8$. It is well-known that
the standar!
 d multiplications of the real algebras of real, complex, quaternionic and
Cayley numbers are both orthogonal multiplications and harmonic morphisms.
Baird also shows that any Clifford system, i.e., an $n$-tuple
$(P_{1},...,P_{n})$ of symmetric endomorphisms of ${\Bbb R}^{2m}$ satisfying $
P_{i}P_{j} + P_{j}P_{i} = 2{\delta}_{ij}Id$ \;for \ $i,j = 1,...n$, defines a
\qh \; $ \varphi:{\Bbb R}^{2m} \longrightarrow{\Bbb R}^{n}$ by
\begin{equation}
\varphi(X) = (\langle P_{1}X,X\rangle,\ldots,\langle
P_{n}X,X\rangle)\label{1c}.
\end{equation}

The author proves \cite{Ou95A} that the ``complete lift'' (Definition 2.1 in
\cite{Ou95A}) of any \qh \;  $ \varphi:{\Bbb R}^{m} \longrightarrow {\Bbb
R}^{n}$ is again a \qh \;$\overline{\varphi}:{\Bbb R}^{2m}\longrightarrow {\Bbb
R}^{n}$, thus provides a method of constructing \qh s from given ones. In
\cite{OuWoo95}, a classification of \qh s is given and it is  also shown that
any umbilical \qh\; arises from a Clifford system. \\

In this work, we solve completely the existence and the classification problem
for \qh s by introducing O-systems of orthogonal transformations of the domain
space.  In $\S 2$ we recall the first classification theorem obtained in
\cite{OuWoo95}, discuss a direct sum construction of harmonic morphisms and
establish a $1-1$ correspondence between equivalence classes of Clifford
systems and umbilical \qh s (Theorem \ref{P1}). In $\S 3$ we introduce the
notion of O-systems and obtain a $1-1$ correspondence between equivalence
classes of Clifford systems and that of O-systems (Theorem \ref{OC}); In $\S 4
$ we establish a $1-1$ correspondence between O-systems and orthogonal
multiplications (Theorem \ref{OO}). Putting these together and using the
Splitting Lemma (Lemma \ref{Split}), we obtain our main theorems (Theorems
\ref{EU} and \ref{EG}). $\S 5$ is devoted to study properties possessed by all
\qh s for fixed pairs of domain and range spaces (including \qh s into ${\Bbb
R}^{2!
 }$, ${\Bbb R}^{3}$, ${\Bbb R}^{4}$, ${\Bbb R}^{5}$, ${\Bbb R}^{8}$, and ${\Bbb
R}^{9}$). Among other things, we show that any \qh \; arises from a single
quadratic function (Theorem \ref{single}). Also we show that we can generalise
the Hopf construction to ``domain-minimal'' but not ``range-maximal'' \qh s
(Theorem \ref{extend}) and finally \qh s in dimensions 2n to n or $n+1$ $(n =
1, 2, 4, 8)$ are shown to be equivalent to standard maps.\\

Applications of O-systems in classifying orthogonal multiplications $F(n,m;m)$
and in constructing isoparametric functions on, and minimal submanifolds of
spheres will appear in the author's paper \cite{Ou95D}.

\section{Quadratic \h s and Clifford systems}

We use $ O({\Bbb R}^{m})$ and $ S({\Bbb R}^{m})$ to denote the set of all
orthogonal endomorphisms and that of all symmetric endomorphisms of ${\Bbb
R}^{m}$ respectively. When the latter is viewed as a Euclidean space it is
understood to have the inner product defined by
\begin{equation}
\langle A,B \rangle = \frac{1}{m} tr(AB) \label{2a}
\end{equation}
for any $A, B \in S({\Bbb R}^{m})$.\\

A {\bf quadratic harmonic morphism}\; $\varphi:{\Bbb R}^{m} \longrightarrow
{\Bbb R}^{n}$ is a \h \; whose components are quadratic functions (i.e.
homogeneous polynomials of degree 2) in $ x_{1},\ldots, x_{m}$. We use
$H_{2}(m,n)$ to denote the set of all \qh s $ \varphi:{\Bbb R}^{m}
\longrightarrow {\Bbb R}^{n}$.\\
Let $ \varphi \in H_{2}(m,n)$ \; with $ ~{\varphi(X) = ( X^{t}A_{1}X,\ldots,
X^{t}A_{n}X)}$. We have proved (\cite{OuWoo95}) that the component matrices
$A_{i}$ have the same rank which is always an even number which we call the
$Q$-rank of $\varphi$. The \qh \; $\varphi$ is said to be {\bf $Q$-nonsingular}
if $Q$-rank($\varphi$) equals the domain dimension, otherwise $\varphi$ is said
to be {\bf $Q$-singular}; We also prove that the $A_{i}$ have the same spectrum
which consists of pairs $\pm \lambda$ of eigenvalues and possibly 0. Then we
obtained the following \\
{\bf Classification Theorem \cite{OuWoo95}:} Let $ \varphi: {\Bbb
R}^{m}\longrightarrow {\Bbb R}^{n} \; (m \geq n)$ be a quadratic harmonic
morphism.\\
(I) If $ \varphi $ is $Q$-nonsingular, then $ m = 2k$ for some $ k \in {\Bbb
N}$ and, with respect to suitable orthogonal coordinates in $ {\Bbb R}^{m}$,
$\varphi$ assumes the normal form
\begin{align}\label{2c}
\varphi(X) =   & \left( X^{t}\left(
                                 \begin{array}{cc}
                                   D  & 0 \\
                                   0 & -D
                                  \end{array}
                              \right)
                                   X,\; X^{t}
                              \left(
\begin{array}{cc}
                                         0  & B_{1} \\
                                        B_{1}^{t} & 0
                                   \end{array}
                             \right) X ,\ldots,  \right. \\
                &   \left.  X^{t} \left(
                                   \begin{array}{cc}
                                     0  & B_{n-1} \\
                                    B_{n-1}^{t} & 0
                                      \end{array}
                                   \right) X \right)
                     \notag
\end{align}
where $B_{i},\; D \in GL({\Bbb R}, k) $ with $D$ diagonal having the positive
eigenvalues as diagonal entries, and satisfy
\begin{equation}\label{2d}
\begin{cases}
DB_{i} = B_{i}D\\
B_{i}^{t}B_{i} = D^{2}\\
B_{i}^{t}B_{j} = - B_{j}^{t}B_{i}.\; \; (i,j,= i,...,n-1, i \neq j).
\end{cases}
\end{equation}
(II) Otherwise $Q$-rank$(\varphi) = 2k$ for some $ k,\  0 \leq k < m / 2$, and
$ \varphi$ is the composition of an orthogonal projection $\pi : {\Bbb
R}^{m}\longrightarrow {\Bbb R}^{2k}$ followed by a $Q$-nonsingular quadratic
harmonic morphism  ${\varphi}_{1}: {\Bbb R}^{2k}\longrightarrow {\Bbb
R}^{n}$.\\

It follows that any quadratic harmonic morphism from an odd-dimen-sional space
is the composition of an orthogonal projection followed by a $Q$-nonsingular
quadratic harmonic morphism from an even-dimensional space. Thus to study
quadratic harmonic morphisms it suffices to consider $Q$-nonsingular ones from
even-dimensional spaces.\\

\begin{definition}\label{D1}
 Let $({M},g)$ and $({N},h)$ be two Riemannian manifolds. Suppose that $
\varphi:M \longrightarrow {\Bbb R}^{n}$ and $\tilde{\varphi}:N \longrightarrow
{\Bbb R}^{n}$ are two $C^{\infty}$ maps. Then the {\bf direct sum} of $
\varphi$ and $\tilde{\varphi}$ is a map
\begin{equation}
\varphi \oplus \tilde{\varphi}: M \times N \longrightarrow {\Bbb R}^{n} \notag
\end{equation}
defined by
\begin{equation}
(\varphi \oplus \tilde{\varphi})(p,q) = \varphi(p) + \tilde{\varphi}(q)\notag
\end{equation}
where $ M \times N$ is the product of $M$ and $N$, endowed with  the Riemannian
product metric $G = (g, h)$.
\end{definition}

For more results on the direct sum construction of \h s see \cite{Ou95C}.

\begin{remark}
It follows from Ou \cite{Ou95C} that the direct sum of any two \h s is again a
\h . In particular, the direct sum of two \qh s $ \varphi:{\Bbb R}^{m}
\longrightarrow {\Bbb R}^{n}$ with $ ~{\varphi(X) = ( X^{t}A_{1}X,\ldots,
X^{t}A_{n}X)}$ \; and \; $\tilde{\varphi}:{\Bbb R}^{l}\longrightarrow {\Bbb
R}^{n}$ \; with \; $\tilde{\varphi}(X) = ( X^{t}B_{1}X,\ldots, X^{t}B_{n}X)$ is
a \qh \\
$\varphi \oplus \tilde{\varphi}:{\Bbb R}^{m+l}\longrightarrow {\Bbb R}^{n}$
given by
\begin{align}\notag
  &(\varphi \oplus \tilde{\varphi})(X,Y) = (X^{t}A_{1}X + Y^{t}B_{1}Y,\ldots,
X^{t}A_{n}X+Y^{t}B_{n}Y)\\
= &\left( (X^{t}\;Y^{t})\left(
                           \begin{array}{cc}  A_{1} & 0 \\ 0 &
B_{1}\end{array}\right)
                           \left(\begin{array}{c} X \\ Y
\end{array}\right),\ldots, (X^{t}\;Y^{t})\left(
                           \begin{array}{cc}  A_{n} & 0 \\ 0 &
B_{n}\end{array}\right)
                           \left(\begin{array}{c} X \\ Y
\end{array}\right)\right)\notag
\end{align}
\end{remark}

A \qh \; is said to be {\bf separable} if it can be written as the direct sum
of two \qh s from smaller dimensional domain spaces. Concerning the existence
of \qh s we note that if there exists $\varphi \in H_{2}(m,n)$ then, by
projection, $H_{2}(m,k) \neq \emptyset $ for $ 1 \leq k <n$; On the other hand,
by using the direct sum construction, we know that if $H_{2}(m,n) \neq
\emptyset$ then $H_{2}(km,n) \neq \emptyset$ for $ k \in {\Bbb N}$.

\begin{definition} \label{D2}

$i)$ A \qh\;$\varphi \in H_{2}(m,n)$ is said to be {\bf range-maximal} if for
fixed $m$, $n$ is the largest range dimension such that $H_{2}(m,n) \neq
\emptyset$ ; It is said to be {\bf domain-minimal} if for fixed $n$, $m$ is the
smallest domain dimension such that $H_{2}(m,n) \neq \emptyset$.\\
$ii)$ Two \qh s $ \varphi, \tilde{\varphi} \in H_{2}(m,n)$ are said to be {\bf
domain-equivalent}, denoted by $\varphi \overset{d}{\sim} \tilde{\varphi}$, if
there exists an isometry $G \in O({\Bbb R}^{m})$ such that $\varphi =
\tilde{\varphi} \circ G$. They are said to be {\bf bi-equivalent}, denoted by
$\varphi \overset{bi}{\sim} \tilde{\varphi}$, if there exist isometries  $G \in
O({\Bbb R}^{m}), H\in O({\Bbb R}^{n})$ such that $\varphi = H^{-1} \circ
\tilde{\varphi} \circ G $.
\end{definition}

Obviously, a \qh \; is domain-minimal if and only if it is not separable. Also,
a domain-minimal \qh \;is always $Q$-nonsingular.

\begin{definition}(see, e.g. \cite{FerKarMun81}) \label{DC}

$i)$ A $(2m, n)$-dimensional {\bf Clifford system} is an $n$-tuple
$(P_{1},\ldots, P_{n})$, denoted by $\{ P_{i}\}$ for short, of symmetric
endomorphisms of ${\Bbb R}^{2m}$ satisfying
\begin{equation}\label{2b}
 P_{i}P_{j} + P_{j}P_{i} = 2{\delta}_{ij}Id \; (i,j = 1,...,n).
\end{equation}
The set of all $(2m,n)$-dimensional Clifford systems is denoted by $ C(2m,n)$.
$ii)$ A {\bf representation} of a Clifford system $\{ P_{i}\} \in C(2m,n)$ is a
$n$-tuple $(A_{1},\ldots, A_{n})$ of symmetric matrices such that, with respect
to some orthonormal basis in ${\Bbb R}^{2m}$, $A_{i}$ is the representation of
$P_{i}$ for $i = 1,...,n$ respectively. A Clifford system is sometimes
specified by its representation as $\{ A_{i} \}$.\\
$iii)$ Let $\{ P_{i} \} \in C(2m,n)$ and $ \{ Q_{i} \} \in C(2l,n)$, then $\{
P_{i}\oplus Q_{i} \}$ is a Clifford system on  ${\Bbb R}^{2m+2l}$, the
so-called {\bf direct sum} of $\{ P_{i} \}$ and $\{ Q_{i} \}$.
\\
$iv)$ A Clifford system $\{ P_{i}\} \in C(2m,n)$ is said to be {\bf irreducible
} if it is not possible to write  ${\Bbb R}^{2m}$ as a direct sum of two
non-trivial subspaces which are invariant under all $P_{i}$.\\
$v)$ Two Clifford systems $\{ P_{i}\},\; \{ Q_{i}\} \in C(2m,n)$ are said to be
{\bf algebraically equivalent}, denoted by $\{ P_{i}\} \overset{a}{\sim} \{
Q_{i}\}$, if there exists $ A \in O({\Bbb R}^{2m})$ such that $ Q_{i} =
AP_{i}A^{t}$ for all $ i = 1,\ldots,n $. They are said to be {\bf geometrically
equivalent}, denoted by $\{ P_{i}\} \overset{g}{\sim} \{ Q_{i}\}$, if there
exists $B\in O(span\{P_{1},\ldots,P_{n}\})\subset O(S({\Bbb R}^{2m}))$ such
that $ \{ B(P_{i}) \} $ and $\{ Q_{i}\}$ are algebraically equivalent.
\end{definition}

A quadratic harmonic morphism is said to be {\bf umbilical} (see
\cite{OuWoo95}) if all the positive eigenvalues of one (and hence all) of its
component matrices are equal. \\

It follows form Baird \cite{Bai83} that any Clifford system $\{P_{i}\} \in
C(2m,n+1)$ gives rise,  by (\ref{1c}), to an umbilical \qh \;with positive
eigenvalue $1$. On the other hand, we know from \cite{OuWoo95} that, up to a
constant factor, any umbilical \qh \;$\varphi \in H_{2}(2m,n+1)$ is
domain-equivalent to one arising from the Clifford system
\begin{equation}\label{2e}
\left( \left( \begin{array}{cc}
I_{m}  & 0\\ 0 & -I_{m}
\end{array} \right),\  \left( \begin{array}{cc}
0  & { \tau_{1}} \\ { \tau_{1}^{t}} & 0
\end{array} \right)
,\ \ldots,  \left( \begin{array}{cc}
0  & { \tau_{n}} \\ { \tau_{n}^{t}} & 0
\end{array} \right) \right)
\end{equation}
where $\tau_{i} \in O({\Bbb R}^{m})$, and satisfy
\begin{equation}
\tau_{i}^{t}\tau_{j} + \tau_{j}^{t}\tau_{i}= 2{\delta}_{ij}Id \; (i,j =
1,...,n).\label{2f}
\end{equation}

\begin{remark}
Note that an umbilical \qh \;with positive eigenvalue $\lambda$ can always be
normalized, by a change of the scalar, to be the one with positive eigenvalue
$1$. Thus the study of umbilical \qh s reduces to the study of umbilical \qh s
with positive eigenvalue $1$.
\end{remark}

Let $H_{2}^{1}(2m,n)$ denote the subset of $H_{2}(2m,n)$ consisting of all
umbilical \qh s with positive eigenvalue $1$. Then, as we have seen from the
above, we have a map
\begin{equation}\label{m1}
 F: C(2m,n) \longrightarrow H_{2}^{1}(2m,n).
\end{equation}
with $F(\{P_{i}\})$ defined by (\ref{1c}).

\begin{theorem}\label{P1}
Let $F$ be the map defined by (\ref{m1}). Then\\
$(i)$ $\{ P_{i}\} \overset{a}{\sim} \{ Q_{i}\} $ if and only if $F \left(\{
P_{i} \}\right) \overset{d}{\sim} F \left(\{ Q_{i} \}\right)$;\\
$(ii)$ $\{ P_{i} \} \overset{g}{\sim} \{ Q_{i} \}$ if and only if $F \left(\{
P_{i} \}\right) \overset{bi}{\sim} F \left(\{ Q_{i} \}\right)$.\\
$(iii)$ $F$ preserves the direct sum operations in the following sense:
\begin{equation}\label{pf1}
 F \left(\{ Q_{i}\oplus R_{i} \}\right) \overset{d}{\sim} F \left(\{ Q_{i}
\}\right) \oplus F \left(\{ R_{i} \}\right),
\end{equation}
and hence  $F$ preserves reducibility in the sense that $\{ P_{i} \}$ is
irreducible if and only if $F \left(\{ P_{i} \} \right) $ is domain-equivalent
to an umbilical \qh \;which is not separable.
\end{theorem}

\begin{proof}
$(i)$ is obviously true. For $(ii)$ we first note that when $S({\Bbb R}^{m})$
is viewed as a Euclidean space with the inner product defined by Equation
(\ref{2a}), then any $\{ P_{i}\} \in C(2m,n)$ becomes an orthonormal set. Thus
span$\{ P_{1},\ldots,P_{n}\}$ can be identified with ${\Bbb R}^{n}$ provided
with the standard inner product, and hence $O(span\{P_{1},\ldots,P_{n}\}) \cong
O({\Bbb R}^{n})$. Now let $B(P_{i}) = a_{i}^{j}P_{j}$. Then it is easy to see
that $B \in O(span\{P_{1},\ldots,P_{n}\})$ if and only if $ (a_{i}^{j}) \in
O({\Bbb R}^{n})$. It is routine to check that $\{B(P_{i})\}$ is indeed a
Clifford system $\in C(2m,n)$. By definition, we have
\begin{align}
F\left( \{B(P_{i})\}\right)(X) &= (\langle B(P_{1})X,X\rangle,\ldots,\langle
B(P_{n})X,X\rangle)\notag\\
                       &=(a_{1}^{j}\langle
P_{j}X,X\rangle,\ldots,a_{n}^{j}\langle   P_{j}X,X\rangle)\notag\\
                       &= H^{-1}\circ F\left(\{P_{i}\}\right)(X)\notag.
\end{align}
where $ H = (a_{i}^{j})^{-1} \in O({\Bbb R}^{n})$. By Definition \ref{DC}, $\{
P_{i}\} \overset{g}{\sim} \{ Q_{i}\}$ if and only if $\{ B(P_{i})\}
\overset{a}{\sim} \{ Q_{i}\}$, which, by $(i)$, is equivalent to $F\left(\{
B(P_{i})\}\right) \overset{d}{\sim} F\left(\{ Q_{i}\}\right)$, i.e., $F\left(\{
Q_{i}\}\right) = H^{-1}\circ F\left(\{P_{i}\}\right)\circ G $ for some $G \in
O({\Bbb R}^{m})$, which means exactly that $F\left(\{P_{i}\}\right)$ and
$F\left(\{ Q_{i}\}\right)$ are bi-equivalent. This ends the proof of $(ii)$.\\

Now we prove $(iii)$, by using $1)$ of Lemma \ref{L3}. it is easy to check that
Equation (\ref{pf1}) holds for any $\{ Q_{i}\} \in C(2k,n)$ and $\{ R_{i}\} \in
C(2l,n)$. Now suppose that $\{ P_{i}\} \in C(2m,n)$ is irreducible, we are to
prove that $F \left(\{ P_{i} \} \right) $ is domain-equivalent to an umbilical
\qh \;which is not separable. Suppose otherwise. Then we would have
\begin{equation}
F \left(\{ P_{i} \} \right)  \overset{d}{\sim} \varphi_{1} \oplus
\varphi_{2},\notag
\end{equation}
where both $\varphi_{1}$ and $\varphi_{2}$ must be umbilical \qh s with
positive eigenvalue $ 1$ since $F \left(\{ P_{i} \} \right)$ is of this kind.
It follows from \cite{OuWoo95} (Theorem 3.3) that there exist $\{ Q_{i}\} \in
C(2k,n)$ and $\{ R_{i}\} \in C(2l,n)$ such that $F \left(\{ Q_{i} \} \right)
\overset{d}{\sim} \varphi_{1}$ and $F \left(\{ R_{i} \} \right)
\overset{d}{\sim} \varphi_{2}$. But then we would have
\begin{equation}
F \left(\{ P_{i} \} \right)  \overset{d}{\sim} \varphi_{1} \oplus \varphi_{2}
 \overset{d}{\sim} F \left(\{ Q_{i} \}\right) \oplus F \left(\{ R_{i}
\}\right)\overset{d}{\sim} F \left(\{ Q_{i}\oplus R_{i} \}\right). \notag
\end{equation}
This means that $\{ P_{i}\} \overset{a}{\sim} \{ Q_{i}\oplus R_{i} \}$, which
is impossible since $\{ P_{i}\}$ is assumed to be irreducible. On the other
hand, it is obviously true that if $F \left(\{ P_{i} \} \right) $ is
domain-equivalent to an unseparable umbilical \qh \;then $\{ P_{i}\}$ is
irreducible. Thus we obtain $(iii)$, which completes the proof of the theorem.
\end{proof}

\section{O-systems and Clifford systems}

\begin{definition}\label{DO}
$i)$ An $(m, n)$-dimensional {\bf O-system} is an $n$-tuple $(\tau_{1},\ldots,
\tau_{n})$, denoted by $\{ \tau_{i}\}$ for short, of orthogonal endomorphisms
of ${\Bbb R}^{m}$ satisfying
\begin{equation}\label{3a}
 \tau_{i}^{t}\tau_{j} + \tau_{j}^{t}\tau_{i} = 2{\delta}_{ij}Id \;\; (i,j =
1,...,n).
\end{equation}
The set of all $(m,n)$-dimensional O-systems is denoted by $ O(m,n)$.\\
$ii)$ A {\bf representation} of an O-system $\{ \tau_{i}\} \in O(m,n)$ is a
$n$-tuple $(a_{1},\ldots, a_{n})$ of orthogonal matrices such that, with
respect to some orthonormal basis in ${\Bbb R}^{m}$, $a_{i}$ is the
representation of $\tau_{i}$ for $i = 1,...,n$ respectively. An O-system is
sometimes specified by its representation as $\{ a_{i} \}$.\\
$iii)$ Let $\{ \rho_{i} \} \in O(m,n)$ and $ \{ \tau_{i} \} \in O(l,n)$, then
$\{ \rho_{i}\oplus \tau_{i} \}$ is an O-system on  ${\Bbb R}^{m+l}$, the
so-called {\bf direct sum} of $\{ \rho_{i} \}$ and $\{ \tau_{i} \}$.
\\
$iv)$ An O-system $\{ \rho_{i}\} \in O(m,n)$ is said to be {\bf irreducible }
if it is not possible to write  ${\Bbb R}^{m}$ as a direct sum of two
non-trivial subspaces which are invariant under all $\rho_{i}$.\\
$v)$ Two O-systems $\{ \rho_{i}\},\; \{ \tau_{i}\} \in O(m,n)$ are said to be
{\bf algebraically equivalent}, denoted by $\{ \rho_{i}\} \overset{a}{\sim} \{
\tau_{i}\}$, if there exists $ \theta \in O({\Bbb R}^{m})$ such that $ \tau_{i}
= \theta \rho_{i}\theta^{t}$ for all $ i = 1,\ldots,n $. They are said to be
{\bf geometrically equivalent}, denoted by $\{ \rho_{i}\} \overset{g}{\sim} \{
\tau_{i}\}$, if $\{ \varrho (\rho_{i}) \} $ and $\{ \tau_{i}\}$ are
algebraically equivalent for some $\varrho \in
O(span\{\rho_{1},\ldots,\rho_{n}\})\subset O(O({\Bbb R}^{m}))$, where $O({\Bbb
R}^{m})$ is provided with the inner product $\langle \tau,\rho \rangle =
\frac{1}{m} tr({\rho}^{t}\tau)$ for any $\tau, \rho \in O({\Bbb R}^{m})$.
\end{definition}

\begin{remark}
We remark that n-tuples $\{ \tau_{i}\}$ of orthogonal endomorphisms satisfying
\begin{equation}\notag
\tau_{i}\tau_{j} + \tau_{j}\tau_{i}= -2{\delta}_{ij}Id\; (i,j = 1,...,n).
\end{equation}
and n-tuples $\{a_{k}\}$ of skew symmetric endomorphisms satisfying Equation
(\ref{3a}) have been used (see e.g., \cite{Hus66, Conl93, OzeTak75,
FerKarMun81}) to study the representations of Clifford algebras. For examples,
we know (see Lemma 24 in \cite{OzeTak75}) that there exists a bijective
correspondence between the set of equivalence classes of $(n-1)$-tuples
$\{a_{k}\}$ of skew symmetric endomorphisms of ${\Bbb R}^{m}$ satisfying
Equation (\ref{3a}) and the set of orthogonal equivalence classes of
representations $\chi (C_{n-1}, \ast)$ of Clifford algebra $C_{n-1}$. On the
other hand, there is a classical result of Radon, Hurwitz and Eckmann's (see
e.g., \cite{Conl93}) saying that there exists $(\sigma(m) -1)$-tuples of {\bf
skew symmetric and orthogonal} endomorphisms of ${\Bbb R}^{m}$ satisfying
Equation (\ref{3a}). In contrast, our results (see Theorem \ref{ES}) claims the
existence of range-maximal $(m,\sigma(m))$-dimensional O-systems which means
that there exist!
  $\sigma(m)$-tuples of {\bf orthogonal} endomorphisms of ${\Bbb R}^{m}$
satisfying Equation (\ref{3a}).
\end{remark}

In comparing with Ozeki and Takeuchi's result, we can use Definition (\ref{DO})
and the relation between O-systems and Clifford systems to establish the
following
\begin{proposition}
There is a bijective correspondence between the set of algebraic equivalence
classes of $O(m,n)$ and that of orthogonal equivalence classes of
representations $\chi (C_{n-1}, \ast)$ of Clifford algebra $C_{n-1}$.
\end{proposition}

\begin{lemma}\label{L3}
$1)$ Let $\{ A_{\alpha}\}$ and $\{ B_{\alpha}\}$ be representations of $\{
P_{\alpha}\}$ and $\{ Q_{\alpha}\} \in C(2m,n+1)$ respectively. Then the direct
sum $\{ P_{\alpha}\oplus Q_{\alpha} \}$ has a representation of the form
\begin{equation}\label{R}
\left\{\left( \begin{array}{cc}
A_{\alpha}  & 0 \\ 0& B_{\alpha}
\end{array} \right)\right\}.
\end{equation}
$2)$ Let $\{a_{i}\}$ and $\{b_{i}\}$ be representations of $\{ \tau_{i}\}$ and
$\{ \rho_{i}\} \in O(m,n)$ respectively. Then the direct sum $\{ \tau_{i}\oplus
\rho_{i} \}$ has a representation of the form
\begin{equation}\label{r}
\left\{\left( \begin{array}{cc}
a_{i}  & 0 \\ 0& b_{i}
\end{array} \right)\right\}.
\end{equation}
$3)$ $\{ P_{\alpha}\} \in C(2m,n+1)$ $($respectively,$\{ \tau_{i}\} \in
O(m,n))$ is reducible if and only if it has a non-trivial representation of the
form $(\ref{R})$ $($respectively, $(\ref{r}))$.
\end{lemma}
\begin{proof}
The proof of the lemma is trivial and is omitted.
\end{proof}

It can be checked that any Clifford system $\{P_{\alpha}\} \in C(2m,n+1)$ is
algebraically equivalent to one given by Equation (\ref{2e}) with an associated
$n$-tuple $\{\tau_{i}\}$ of orthogonal endomorphisms satisfying (\ref{2f}).
Thus every Clifford system $\{P_{\alpha}\} \in C(2m,n+1)$ corresponds to an
O-system $ \{\tau_{i}\} \in O(m,n)$, and hence we have a surjective map
\begin{equation} \label{3b}
f:C(2m,n+1) \longrightarrow O(m,n),
\end{equation}
with $f\left(\{P_{\alpha}\}\right) = \{\tau_{i}\}$.
\begin{theorem}\label{OC}
Let $f$ be the map defined by (\ref{3b}). Then \\
$(1)$ $f$ induces a bijective correspondence, ${\overline
f}:C(2m,n+1)/\overset{a}{\sim} \longrightarrow O(m,n)/\overset{a}{\sim}$.\\
$(2)$ $f$ preserves the direct sum operations in the following sense:
\begin{equation}\label{psum}
 f \left(\{ Q_{\alpha}\oplus R_{\alpha} \}\right) = f \left(\{ Q_{\alpha}
\}\right) \oplus f \left(\{ R_{\alpha} \}\right),
\end{equation}
and hence  $f$ preserves reducibility in the sense that $\{ P_{\alpha} \}$ is
irreducible if and only if $f \left(\{ P_{\alpha} \} \right) $ is irreducible.
\end{theorem}

\begin{proof}
For $(1)$, it is evident that $f$ induces an onto map
\begin{equation}\notag
 {\overline f}:C(2m,n+1)/\overset{a}{\sim} \longrightarrow
O(m,n)/\overset{a}{\sim}.
\end{equation}
It remains to prove that ${\overline f}$ is $1-1$. To this end, suppose that
\begin{equation}
 {\overline f}([\{ P_{\alpha}\}]) =[\{\tau_{i}\}] = {\overline f}([\{
Q_{\alpha}\}]) =[\{\rho_{i}\}] \notag
\end{equation}
Then we have $\{\tau_{i}\} \overset{a}{\sim} \{ \rho_{i}\}$. Thus there exists
$ \theta \in O({\Bbb R}^{m})$ such that $ \tau_{i} = \theta \rho_{i}\theta^{t}$
for all $ i = 1,\ldots,n $. Now one can check that
\begin{equation}
 A = \left( \begin{array}{cc}
\theta & 0 \\ 0 & \theta
\end{array} \right) \in O({\Bbb R}^{2m})\notag
\end{equation}
and that
\begin{equation}
A\left( \begin{array}{cc}
0  & { \tau_{i}} \\ { \tau_{i}^{t}} & 0
\end{array} \right)A^{t} = \left( \begin{array}{cc}
0  & { \rho_{i}} \\ { \rho_{i}^{t}} & 0
\end{array} \right),\notag
\end{equation}
which means that $\{ P_{\alpha}\}$ and $\{Q_{\alpha}\}$ are algebraically
equivalent, and hence $[\{ P_{\alpha}\}] = [\{ Q_{\alpha}\}]$.\\
For the first statement of $(2)$, we can check that
\begin{equation}\notag
\{Q_{\alpha}\oplus R_{\alpha}\}\overset{a}{\sim} \left\{
\left(\begin{array}{cccc}
I_{m} & 0 & 0 & 0\\
0 & I_{l} & 0 & 0\\
0 & 0 & -I_{m} & 0 \\
0 & 0 & 0 & -I_{l} \end{array}\right),\; \left(\begin{array}{cccc}
0 & 0 & \tau_{i} & 0\\
0 & 0  & 0 & \rho_{i}\\
\tau_{i}^{t} & 0 & 0 & 0 \\
0 & \rho_{i}^{t}& 0 & 0\end{array}\right)\right\}.
\end{equation}
Thus, by definition of $f$ and $(1)$, we obtain Equation (\ref{psum}), which
together with $(1)$ gives the second statement of $(2)$. This completes the
proof of the theorem.
\end{proof}

\begin{remark}\label{3.4}
Combining $(i)$ of Theorem \ref{P1} and $(1)$ of Theorem \ref{OC}, we have the
following $1-1$ correspondences:
\begin{equation}
O(m,n)/\overset{a}{\sim} \overset{1-1}{\longrightarrow}
C(2m,n+1)/\overset{a}{\sim} \overset{1-1}{\longrightarrow}  H_{2}^{1}(2m,n+1)/
\overset{d}{\sim}.\notag
\end{equation}
\end{remark}

\begin{proposition}\label{P3}
Let $\{ \tau_{i}\} \in O(m,n)$ be an O-system. Then\\
$(a)$ $\{ \tau_{i}^{t}\}$ is also an O-system in $O(m,n)$.\\
$(b)$ any subset consisting of $k$ elements of $\{ \tau_{i}\}$ forms an
$(m,k)$-dimensional O-system.
\end{proposition}
\begin{proof}
The proof is a straightforward checking of the defining Equation (\ref{3a}) and
is omitted.
\end{proof}
\begin{remark}\label{R3}
It follows from Proposition \ref{P3} that if $O(m,n) = \emptyset$ then
$O(m,n+p) = \emptyset$ for $ p \geq 1$. On the other hand, by the direct sum
operation, if $O(m,n) \neq \emptyset$ then $O(km,n) \neq \emptyset$ for $ k
\geq 2$. Therefore it is meaningful to put the qualifiers {\bf
``range-maximal''} and {\bf ``domain-minimal''} before  a Clifford system and
an O-system in a similar way as they are used in Definition \ref{D2}.
\end{remark}
\begin{proposition}\label{P4}
For $ n \geq 2, k \in {\Bbb N},\ O(2k+1,n) = \emptyset$. Geometrically, there
exists no $Q$-nonsingular umbilical \qh s $ \varphi: {\Bbb
R}^{4k+2}\longrightarrow {\Bbb R}^{n+1}$ for $ n \geq 2$.
\end{proposition}
\begin{proof}
By Remark \ref{R3}, we need only to show that $O(2k+1,2) = \emptyset$. In fact,
if there were $\{\tau_{1},\tau_{2}\}\in O(2k+1,2)$, then we would have $
\tau_{1},\tau_{2} \in O({\Bbb R}^{2k+1})$ satisfying $\tau_{1}^{t}\tau_{2}= -
\tau_{2}^{t}\tau_{1}$. It follows that
\begin{equation}
\det\tau_{1}\det\tau_{2} = (-1)^{2k+1}\det\tau_{2}\det\tau_{1} = -
\det\tau_{1}\det\tau_{2},\notag
\end{equation}
and hence $\det\tau_{1} =0$ or $\det\tau_{2}= 0$ which contradict the fact that
$\tau_{1},\tau_{2}$ are orthogonal.
\end{proof}

\section{O-systems and Orthogonal multiplications}

\begin{definition}\label{OM}
An {\bf orthogonal multiplication} is an ${\Bbb R}$-bilinear map $\mu :{\Bbb
R}^{p} \times {\Bbb R}^{q} \longrightarrow {\Bbb R}^{n} $ with $\|\mu(x,y)\| =
\|x\|\|y\|$ for all $x \in {\Bbb R}^{p}, y \in {\Bbb R}^{q}$. We use $F(p,q;n)$
to denote the set of all orthogonal multiplications $\mu :{\Bbb R}^{p} \times
{\Bbb R}^{q} \longrightarrow {\Bbb R}^{n}$.
\end{definition}
Though the existence of \om s is a purely algebraic problem it is closely
related to the existence of interesting geometric objects such as vector fields
on spheres and harmonic maps into spheres:\\
{\bf Fact A:}(see \cite{Hus66}) If there exists an \om \; $\mu \in F(n,m;m)$
then there exist $(n-1)$ linearly independent vector fields on ${\Bbb
S}^{m-1}$.\\
{\bf Fact B:}(see \cite{EelLem83}) For any orthogonal multiplication $\mu \in
F(p,q;n)$, the restriction of $\mu$ provides a harmonic map ${\Bbb S}^{p-1}
\times {\Bbb S}^{q-1} \longrightarrow {\Bbb S}^{n-1}$. If $ p = q$, the Hopf
construction
\begin{equation}\label{4a}
H(x,y) = \left(\|x\|^{2} - \|y\|^{2},\  2\mu(x,y)\right)\; :{\Bbb R}^{p} \times
{\Bbb R}^{p} \longrightarrow {\Bbb R}^{n+1}.
\end{equation}
provides a harmonic map ${\Bbb S}^{2p-1}\longrightarrow {\Bbb S}^{n}$.\\

In relating to \qh s, Baird has proved (\cite{Bai83} Theorem 7.2.7) that an \om
\;$\mu \in F(p,q;n)$ is a \h \ if and only if the dimensions $ p= q =n$ and $n
= 1,2, 4$, or $8$.\\

Now we give a link between O-systems and \om s as
\begin{theorem}\label{OO}
There exists a $1-1$ correspondence between $O(m,n)$ and $F(n,m;m)$. The
correspondence is given by
\begin{equation}
O(m,n)\ni \{\tau_{i}\}\longmapsto \mu_{\tau} :{\Bbb R}^{n} \times {\Bbb R}^{m}
\longrightarrow {\Bbb R}^{m},\notag
\end{equation}
with
\begin{equation}\label{4b}
\mu_{\tau}(x,y) = x^{i}\mu_{\tau}(e_{i},y) = x^{i}\tau_{i}(y).
\end{equation}
where $x = x^{i}e_{i}$ and $\{e_{i}\}$ is the standard basis for ${\Bbb
R}^{n}$.
\end{theorem}

\begin{proof}
It is an elementary fact from linear algebra that there is a $1 - 1$
correspondence, given by (\ref{4b}), between the set of $n$-tuples $\{
\tau_{i}\}$ of linear endomorphisms of ${\Bbb R}^{m}$ and that of bilinear maps
$\mu :{\Bbb R}^{n} \times {\Bbb R}^{m} \longrightarrow {\Bbb R}^{m}$. It
remains to show that $\{ \tau_{i}\}$ is an O-system if and only if the
corresponding $\mu_{\tau}$ is an \om . To this end, we first note the following
two facts:\\
$(I)$ A linear transformation $\tau :{\Bbb R}^{m} \longrightarrow {\Bbb R}^{m}$
is orthogonal if and only if $\tau$ sends ${\Bbb S}^{m-1}$ into ${\Bbb
S}^{m-1}$;\\
$(II)$ A bilinear map $\mu :{\Bbb R}^{n} \times {\Bbb R}^{m} \longrightarrow
{\Bbb R}^{m} $ is an \om \; if and only if $\mu $ sends ${\Bbb S}^{n-1} \times
{\Bbb S}^{m-1}$ into ${\Bbb S}^{m-1}$.\\
Now suppose that $\{ \tau_{i}\} \in O(m,n)$. Then $x^{i}\tau_{i}$ is orthogonal
for any $(x^{1},\ldots,x^{n}) \in {\Bbb S}^{n-1}$ since
\begin{equation}\label{4c}
(x^{i}\tau_{i})(x^{i}\tau_{i})^{t} = \sum_{i=1}^{n} (x^{i})^{2}Id + \sum_{i
\neq j} x^{i}x^{j}(\tau_{i}^{t}\tau_{j} + \tau_{j}^{t}\tau_{i})= Id.
\end{equation}
Therefore, for any $x^{i}e_{i} \in {\Bbb S}^{n-1}$ and $y \in {\Bbb S}^{m-1}$
we have

\begin{equation}\label{4d}
\|\mu_{\tau}(x,y)\| = \|x^{i}\tau_{i}(y)\| = \|(x^{i}\tau_{i})(y)\| = 1.
\end{equation}
{}From this and $(II)$ it follows that $\mu_{\tau}$ is an \om . Conversely, if
$\mu_{\tau}$ is an \om \;then by $(II)$, Equation (\ref{4d}) holds for any
$x^{i}e_{i} \in {\Bbb S}^{n-1}$ and $y \in {\Bbb S}^{m-1}$. Thus, by $(I)$,
$x^{i}\tau_{i}$ is orthogonal for any $(x^{1},\ldots,x^{n}) \in {\Bbb
S}^{n-1}$. It follows that Equation(\ref{4c}) holds for arbitrary
$(x^{1},\ldots,x^{n}) \in {\Bbb S}^{n-1}$, which implies that
$\tau_{i}^{t}\tau_{j} + \tau_{j}^{t}\tau_{i} = 0$ for $ i,j = 1,\ldots, n$ and
$i\neq j$. Thus $\{ \tau_{i}\} \in O(m,n)$, which ends the proof of the
theorem.
\end{proof}

\begin{remark}
It is known (Smith \cite{Smi72}, see also Eells and Ratto \cite{EelRat93}) that
there is a $1-1$ correspondence between $F(p,q;n)$ and the set of totally
geodesic embeddings of\; ${\Bbb S}^{p-1}$ into \;$O_{n,q}$, the Stiefel
manifold of orthogonal q-frames in n-space with suitable normalization. In
particular, there is a $1-1$ correspondence between $F(n,m;m)$ and geodesic
$(n-1)$-spheres in $O({\Bbb R}^{m})$.
\end{remark}

Now we are ready to give the following existence theorem for Clifford systems
and O-systems.
\begin{theorem}\label{ES}
$A)$ Let $\sigma(m) = 2^{c} + 8d$ for any $m\in{\Bbb N}$, uniquely written as
$m = (2r + 1)2^{c+4d} \ (r, d \geq 0,\ 1\leq c \leq 3)$. Then there exist
range-maximal $(m,\sigma(m))-$dimensional O-systems, and hence range-maximal
$(2m,\sigma(m)+1)$-dimensional Clifford systems.\\
$B)$ For any $ n \in {\Bbb N}$ there exist domain-minimal
$(2m(n),n+1)$-dimensional Clifford systems, and hence domain-minimal
$(m(n),n)$-dimensional O-systems for and only for the  $(m(n),n)$ values listed
in Table 1.
\begin{table}[htb]
\centering
\begin{tabular}{|c|c|c|c|c|c|c|c|c|c|c|}
\hline n & 1 & 2 & 3 & 4 & 5 & 6 & 7 & 8 & \ldots & n+8 \\ \hline
m(n) & 1 & 2 & 4 & 4 & 8 & 8 & 8 & 8 & \ldots & 16 m(n) \\ \hline
\end{tabular}
\bigskip
\caption{}
\end{table}
\end{theorem}

\begin{proof}
{}From Theorem \ref{OO} we know that there exists an $(m,n)$-dimen\-sional
O-system if and only if there exists an \om \;$\mu \in F(n,m;m)$. Now the
existence of range-maximal O-systems follows from a classical result of Hurwitz
\cite{Hurw23} (see also Radon \cite{Rado22} and Eckmann \cite{Eckm52}) that for
any $m \in {\Bbb N} $, uniquely written as $m = (2r + 1)2^{c+4d} \ (r, d \geq
0,\ 1\leq c \leq 3)$, there exist \om s $\mu \in F(n,m;m)$ for $n = \sigma (m)
= 2^{c}+8d$ which is also the largest number possible for such \om s to exist.
The existence of range-maximal Clifford systems then follows from Theorem
\ref{OC}. Thus we obtain $A)$.\\
For $B)$, we note that a Clifford system or an O-system is domain-minimal if
and only if it is irreducible. It follows from \cite{FerKarMun81} that
$(2m(n),n+1)$-dimensional irreducible Clifford systems exist precisely for the
values of $(m(n),n)$ listed in Table 1. Again the relation between Clifford
systems and O-systems (Theorem \ref{OC}) gives the existence of domain-minimal
O-systems, which completes the proof of the theorem.
\end{proof}

\begin{theorem}{$($Existence of umbilical quadratic harmonic
morphisms$)$}\label{EU}
$(a)$ Let $\sigma(m) = 2^{c} + 8d$ for any $m\in{\Bbb N}$, uniquely written as
$m = (2r + 1)2^{c+4d} \ (r, d \geq 0,\ 1\leq c \leq 3)$. Then there exist
range-maximal $Q$-nonsingular umbilical \qh s $ \varphi: {\Bbb
R}^{2m}\longrightarrow {\Bbb R}^{\sigma (m)+1}$.\\
$(b)$ For any $ n \in {\Bbb N}$ there exist domain-minimal $Q$-nonsingular
umbilical \qh s $ \varphi: {\Bbb R}^{2m(n)}\longrightarrow {\Bbb R}^{n+1}$ for
and only for the  $(m(n),n)$ values listed in Table 1.\\
$(c)$ Any other $Q$-nonsingular umbilical \qh s into ${\Bbb R}^{n+1}$ exist
precisely in the cases
\begin{equation}
{\Bbb R}^{2km(n)}\longrightarrow {\Bbb R}^{n+1},\; k \geq 2,\notag
\end{equation}
where they are domain-equivalent to a direct sum of some domain-minimal
umbilical \qh \;\ ~{$ \varphi: {\Bbb R}^{2m(n)}\longrightarrow {\Bbb
R}^{n+1}$}.
\end{theorem}

\begin{proof}
Using Theorem \ref{ES} and the map (\ref{m1}) we obtain  statements $(a)$ and
$(b)$ immediately. For $(c)$ we first note from \cite{OuWoo95} that any
$Q$-nonsingular umbilical \qh \;$\varphi$ is domain-equivalent to $\lambda
\varphi_{0}$ for $\varphi_{0}$ given by a Clifford system $\{P_{i}\} \in
C(2km(n),n+1)$. It is known (see e.g. \cite{FerKarMun81}) that any Clifford
system is algebraically equivalent to a direct sum of irreducible ones. Thus,
according to Table 1, any Clifford system $\{P_{i}\} \in C(2km(n),n+1)$ is
algebraically equivalent to
\begin{equation}
\{ P_{i}^{1} \oplus \ldots \oplus P_{i}^{k} \},\notag
\end{equation}
where $\{ P_{i}^{\alpha}\}\in C(2m(n),n+1)$ is irreducible for all $\alpha =
1,\ldots,k$. By using Theorem \ref{P1} we see that
\begin{equation}
\varphi_{0} \overset{d}{\sim} F \left(\{ P_{i}^{1} \}\right) \oplus \ldots
\oplus F \left(\{ P_{i}^{k} \}\right),\notag
\end{equation}
which gives $(c)$ and hence we obtain the theorem.
\end{proof}

As an immediate consequence, we have

\begin{corollary}\label{4.5}
Let $\varphi: {\Bbb R}^{2m}\longrightarrow {\Bbb R}^{n}$ be a $Q$-nonsingular
umbilical \qh\;. Then either\\
 $(a)$ $\varphi$ is domain-minimal, in which case, $ \varphi
\overset{d}{\sim}\lambda \varphi_{0}$ for $\varphi_{0} \in H_{2}^{1}(2m,n)$
given by an irreducible Clifford system, or\\
 $(b)$ $\varphi \overset{d}{\sim} \lambda (\varphi_{1}\oplus \ldots \oplus
\varphi_{k})$, where all $\varphi_{i} \in H_{2}^{1}(2l,n) \;(kl = m)$ are
domain-minimal.
\end{corollary}

Now we give a splitting lemma which will give the existence of general \qh s.
\begin{lemma}{\bf (The Splitting Lemma)}\label{Split}
Let $ \varphi: {\Bbb R}^{2m}\longrightarrow {\Bbb R}^{n}$ be a $Q$-nonsingular
\qh\;. Then either\\
$(i)$ $\varphi$ is umbilical, or\\
$(ii)$ $\varphi$ is domain-equivalent to a direct sum of $Q$-nonsingular
umbilical \qh s.
\end{lemma}
\begin{proof}
Let $ \varphi: {\Bbb R}^{2m}\longrightarrow {\Bbb R}^{n}$ be a $Q$-nonsingular
\qh\;. We will prove the lemma by checking the following two cases:\\
{\bf Case I: All the positive eigenvalues $ \lambda_{1},\ldots,\lambda_{m}$ are
distinct.}\\

{\bf Claim 1:} For $n\geq 3$, there exists no $Q$-nonsingular \qh \; $\varphi:
{\Bbb R}^{2m}\longrightarrow {\Bbb R}^{n}$ with all positive eigenvalues
distinct.\\
{\bf Proof of Claim 1:} Since the composition of a \qh \;with distinct positive
eigenvalues followed by a projection is again a \qh \;of this kind, it is
enough to do the proof for the case ${\Bbb R}^{2m}\longrightarrow {\Bbb
R}^{3}$. Suppose otherwise, if there were a \qh \;of this kind, then by the
Classification Theorem, we would have
\begin{equation}\label{4e}
D = \left(
\begin{array}{cccc}
\lambda_{1} & 0 & \ldots & 0 \\
0 &\lambda_{2} & \ldots & 0 \\
\vdots & \vdots & \ddots & \vdots \\
0 & 0 & \ldots & \lambda_{m}
\end{array}
\right),\;\; B_{1} = \left( b_{ij}^{1} \right), \; B_{2} = \left(
b_{ij}^{2}\right) \in GL({\Bbb R},m)
\end{equation}
and satisfy Equation (\ref{2d}). Now from the first equation of (\ref{2d}) we
know that $B_{i}$ must be of diagonal form
\begin{equation}
B_{i} = \left(\begin{array}{cccc}
b_{11}^{i} & 0 & \ldots & 0 \\
0 & b_{22}^{i} & \ldots & 0 \\
\vdots & \vdots & \ddots & \vdots \\
0 & 0 & \ldots & b_{mm}^{i}
\end{array}
\right).\notag
\end{equation}
But then the third equation of (\ref{2d}) says that $ b_{jj}^{1}b_{jj}^{2} = 0$
for $j = 1,\ldots,m$. This implies that either $b_{jj}^{1} =0$ or $b_{jj}^{2}
=0$, which is impossible since $B_{i} \;(i=1,2)$ is  non-singular. This ends
the proof of Claim 1..

{\bf Claim 2:} Any $Q$-nonsingular \qh \; \\$\varphi: {\Bbb
R}^{2m}\longrightarrow {\Bbb R}^{2}$ with all positive eigenvalues distinct is
domain-equivalent to a direct sum of umbilical ones.\\
{\bf Proof of Claim 2:} We know from \cite{OuWoo95} that any $Q$-nonsingular
\qh \; $\varphi: {\Bbb R}^{2m}\longrightarrow {\Bbb R}^{2}$ is
domain-equivalent to the normal form
\begin{equation}\label{4f}
\tilde{\varphi}(X) = \left( X^{t}\left(
                                 \begin{array}{cc}
                                   D  & 0 \\
                                   0 & -D
                                  \end{array}
                              \right)
                                   X,\; X^{t}
                              \left(\begin{array}{cc}
                                         0  & B_{1} \\
                                        B_{1}^{t} & 0
                                   \end{array}
                             \right) X\right),
\end{equation}
where $D$ and $B_{1}$ are as in (\ref{4e}) and satisfy Equation (\ref{2d}).
{}From Equation (\ref{2d}) and the hypothesis on $\lambda_{i}'s$ we deduce that
$B_{1}$ must be diagonal form
\begin{equation}
B_{1} = \left(
\begin{array}{cccc}
\pm \lambda_{1} & 0 & \ldots & 0 \\
0 & \pm \lambda_{2} & \ldots & 0 \\
\vdots & \vdots & \ddots & \vdots \\
0 & 0 & \ldots & \pm \lambda_{m}
\end{array}
\right).\notag
\end{equation}
Inserting this into (\ref{4f}), we have
\begin{align}
\tilde{\varphi}(X) = ( & \lambda_{1}(x_{1}^{2} - x_{m+1}^{2})+ \ldots +
\lambda_{m}(x_{m}^{2} - x_{2m}^{2}),\;\notag\\
 & \pm 2\lambda_{1}x_{1}x_{m+1} \pm \ldots \pm 2\lambda_{m}x_{m}x_{2m} ).\notag
\end{align}
It is easy to check that
\begin{equation}
\varphi \overset{d}{\sim} \tilde{\varphi} \overset{d}{\sim}
\lambda_{1}\varphi_{0}\oplus \ldots \oplus \lambda_{m}\varphi_{0},\notag
\end{equation}
where $\varphi_{0}:{\Bbb R}^{2}\cong {\Bbb C}\longrightarrow {\Bbb R}^{2}\cong
{\Bbb C}$ with $\varphi_{0}(z) = z^{2}$ which is clearly umbilical. Thus we
obtain Claim 2. Combining Claim 1 and 2 we see that the lemma is true for Case
I.

{\bf Case II: $\varphi$ has some equal positive eigenvalues.}\\

Without loss of generality, we may assume that the positive eigenvalues satisfy
$ \lambda_{1} = \ldots = \lambda_{k}\neq \lambda_{l}\;(k < l \leq m)$. By the
Classification Theorem, we know that $\varphi\overset{d}{\sim} \tilde{\varphi}$
given by the normal form (\ref{2c}) with $B_{i}\;(i=1, \ldots, n-1),\;D \in
GL({\Bbb R},m)$ satisfying (\ref{2d}). Now using Equation (\ref{2d}) and the
hypothesis on $\lambda_{i}$'s we can check that $B_{i}$ must take the form
\begin{equation}
B_{i}= \left(\begin{array}{cc}
b_{i} & 0\\
0 & c_{i}\end{array} \right),\; i = 1, \ldots, n-1,\notag
\end{equation}
where $b_{i} \in GL({\Bbb R},k),\;c_{i} \in GL({\Bbb R},m-k)$. Writing
\\$(\underbrace{x_{1},\ldots, x_{k}},\;
\underbrace{x_{k+1},\ldots,x_{m}},\;\underbrace{x_{m+1},\ldots,x_{m+k}},\;\underbrace{x_{m+k+1},\ldots,x_{2m}})$ \\as $(X_{1}, X_{2}, X_{3}, X_{4})$, we can check that after an orthogonal change of the coordinates of the form

\begin{equation}
\left(\begin{array}{c}
X_{1}\\X_{2}\\X_{3}\\X_{4}\end{array}\right) =\left(\begin{array}{cccc}
I_{k} & 0 & 0 & 0\\
0 & 0  & I_{m-k} & 0\\
0 & I_{k} & 0 & 0 \\
0 & 0 & 0 & I_{m-k} \end{array}\right) \left(\begin{array}{c}
\tilde{X_{1}}\\ \tilde{X_{2}}\\ \tilde{X_{3}}\\
\tilde{X_{4}}\end{array}\right),\notag
\end{equation}
we have, $\tilde{\varphi}(\tilde{X}) = \varphi_{1}(\tilde{X_{1}},
\tilde{X_{2}}) + \overline{\varphi}_{2}(\tilde{X_{3}}, \tilde{X_{4}})$, i.e.,
$\tilde{\varphi}\overset{d}{\sim}\varphi_{1} \oplus \overline{\varphi}_{2}$
with $\varphi_{1}: {\Bbb R}^{2k}\longrightarrow {\Bbb R}^{n}$ and
$\overline{\varphi}_{2}:{\Bbb R}^{2(m-k)}\longrightarrow {\Bbb R}^{n}$ both
quadratic maps. By using Theorem 1.5 in \cite{OuWoo95} and the fact that
$\tilde{\varphi}(\tilde{X})$ is a Q-nonsingular \qh \;we can prove that both
$\varphi_{1}$ and $\overline{\varphi}_{2}$ are Q-nonsingular \qh s with
$\varphi_{1}$ umbilical since $\lambda_{1},\ldots,\lambda_{k}$ are supposed to
be equal. Now the same process applies to $\overline{\varphi}_{2}$ and so on
until we obtain
\begin{equation}
\varphi \overset{d}{\sim} \tilde{\varphi} \overset{d}{\sim} \varphi_{1}\oplus
\ldots \oplus \varphi_{i},\notag
\end{equation}
for all umbilical $\varphi_{1}, \ldots ,\varphi_{i}$.\\
Therefore, we have seen from Case I and Case II that either $\varphi$ is
umbilical or, $\varphi$ is domain-equivalent to a direct sum of $Q$-nonsingular
umbilical ones. This ends the proof of the lemma.
\end{proof}

Recall that a domain-minimal \qh \;is always $Q$-nonsingular and not separable,
we see immediately, from the Splitting Lemma, the following

\begin{corollary}\label{4.7}
Any domain-minimal \qh \; \\$\varphi: {\Bbb R}^{2m}\longrightarrow {\Bbb
R}^{n}$ is umbilical.
\end{corollary}

{}From the Splitting Lemma and Theorem \ref{EU} we obtain the following
existence theorem for general \qh s.
\begin{theorem}\label{EG}
For $n \in{\Bbb N}$, $Q$-nonsingular \qh s exist precisely in the cases ${\Bbb
R}^{2km(n)}\longrightarrow {\Bbb R}^{n+1}$ for $k \in{\Bbb N}$, and $m(n)$
depending on $n$ given by Table 1. Furthermore\\
$(I)$ If $k = 1$, then $\varphi$ is domain-minimal and hence umbilical, and
$\varphi \overset{d}{\sim} \lambda \varphi_{0}$ for $\varphi_{0} \in
H_{2}^{1}(2m(n),n+1)$.  Otherwise\\
$(II)$ $\varphi$ is domain-equivalent to a direct sum of $k$ domain-minimal \qh
s, i.e.,
\begin{equation}
\varphi \overset{d}{\sim} \lambda_{1}\varphi_{1}\oplus \ldots \oplus
\lambda_{k}\varphi_{k},\label{4g}
\end{equation}
for $\varphi_{i} \in H_{2}^{1}(2m(n),n+1)$, given by irreducible Clifford
systems.
\end{theorem}

\begin{remark}\label{xx}
Note that $\lambda_{1},\ldots, \lambda_{k}$ in (\ref{4g}) are the distinct
positive eigenvalues of $\varphi$. If we allow some but not all of them to be
zero then our results include also $Q$-singular cases. Thus Equation (\ref{4g})
gives, up to domain-equivalence, a general form of \qh s ${\Bbb
R}^{2km(n)}\longrightarrow {\Bbb R}^{n+1}$.
\end{remark}

\begin{corollary}
Let $\varphi:{\Bbb R}^{2km(n)}\longrightarrow {\Bbb R}^{n+1}$ be a \qh \;.
Then\\
$i)$ for $ n \not\equiv 0 \; \mod \;4$, $\varphi
\overset{d}{\sim}\lambda_{1}\varphi_{0}\oplus \ldots \oplus
\lambda_{k}\varphi_{0}$, where $\varphi_{0}:{\Bbb R}^{2m(n)}\longrightarrow
{\Bbb R}^{n+1}$ is a domain-minimal \qh \;given by an irreducible Clifford
system. Otherwise\\
$ii)$ $ n \equiv 0 \; \mod \;4$, $\varphi$ belongs to one of $2^{k-1}$
bi-equivalent classes.
\end{corollary}

\begin{proof}
{}From \cite{FerKarMun81}, we know that there exists only one algebraic
equivalent class in $C(2m(n),n+1)$ for $ n \not\equiv 0 \; \mod \;4$, and two
for $ n \equiv 0 \; \mod \;4$. The the statement $i)$ now follows from Theorem
\ref{EG} and Remark \ref{xx}.. For statement $ii)$, we first note (see
\cite{FerKarMun81}) that two algebraically different irreducible Clifford
systems differ only by a minus sign before one, say the last, of their
elements. Correspondingly, two domain-minimal \qh s\;$\lambda \varphi_{1}$ and
$\lambda \varphi_{2}$ from two different domain-equivalent classes differ  only
by a minus sign before, say,  the last component functions. Therefore, in
constructing direct sum of the form $\lambda_{1}\varphi_{1}\oplus \ldots \oplus
\lambda_{k}\varphi_{k}$ with fixed tuple $(\lambda_{1}, \ldots,\lambda_{k})$,
we get $2^{k}$ possibilities, of which, half can be obtained from the other by
an orthogonal change of the coordinates in the range space as one can check.!
  Thus $ii)$ follows, which completes the proof of the corollary.
\end{proof}

\section{Properties of \qh s}

\begin{definition}
Let $(M, g)$ be a space form (i.e., either Euclidean sphere ${\Bbb S}^{m}$,
Euclidean space ${\Bbb R}^{m}$ or hyperbolic space ${\Bbb H}^{m}$). A smooth
function $f :M \longrightarrow {\Bbb R}$ is called {\bf isoparametric} if
\begin{align}
& \|d\, f(x)\|^{2} = \psi_{1}(f(x)),\notag\\
& \triangle f(x) = \psi_{2}(f(x)),\notag
\end{align}
for some smooth functions $\psi_{1},\; \psi_{2}:{\Bbb R} \longrightarrow {\Bbb
R}$.
\end{definition}

Such functions were introduced by Cartan \cite{Car38} in 1938. Their
description on Euclidean space and hyperbolic space is relatively trivial, but
on the sphere they are rich in geometry. More recent studies on such functions
have been made in \cite{FerKarMun81, Mun80, Nom75, OzeTak75, OzeTak76,
Takagi76, TakTak72}.

It is well-known that all isoparametric functions $f :{\Bbb
S}^{m-1}\longrightarrow {\Bbb R}$ arise from the restriction of a homogeneous
polynomial $F:{\Bbb R}^{m}\longrightarrow {\Bbb R}$ of degree $p$ with
\begin{align}
&\|\bigtriangledown F\|^{2} = p^{2}\|x\|^{2p-2},\label{5a}\\
&\triangle F = c \|x\|^{p-2}.\label{5b}
\end{align}
where $c = 0 $ if the multiplicities of the distinct principal curvatures are
equal.

Suppose we are given a homogeneous polynomial $F:{\Bbb R}^{m}\longrightarrow
{\Bbb R}$ of degree $p$ with $c = 0$, satisfying (\ref{5a}) and (\ref{5b}).
Given any $A \in O({\Bbb R}^{m})$, we define $G:{\Bbb R}^{m}\longrightarrow
{\Bbb R}$ by putting $G = F \circ A$. Then $G$ is also a polynomial satisfying
(\ref{5a}) and (\ref{5b}). Baird noted that if $A \in O({\Bbb R}^{m})$ can be
so chosen that $ \langle \bigtriangledown G_{X}, \bigtriangledown F_{X} \rangle
= 0 $ for any $X \in {\Bbb R}^{m}$, then $\varphi = (F, G) $ gives a nontrivial
\h \;defined by homogeneous polynomials of degree $p$. This provides a possible
way to construct polynomial \h s from a single homogeneous polynomial. However,
this method fails in a case of homogeneous degree $4$ polynomial \h \;${\Bbb
R}^{6}\longrightarrow {\Bbb R}^{2}$ (see Theorem 8.3.5 in \cite{Bai83}). We
shall prove that this method works for any \qh s.

\begin{theorem}\label{single}
Any \qh \;$\varphi:{\Bbb R}^{2km(n)}\longrightarrow {\Bbb R}^{n+1}$ arises from
a single quadratic function $F:{\Bbb R}^{2km(n)}\longrightarrow {\Bbb R}$,\; $
F = \lambda_{1} F_{0} \oplus \ldots \oplus \lambda_{k}F_{0}$, where $
\lambda_{i} \geq 0$ are constants with at least one not zero, and $F_{0}:{\Bbb
R}^{2m(n)}\longrightarrow {\Bbb R}$,
\begin{equation}
F_{0}(x_{1}, \ldots, x_{2m}) = {x_{1}}^{2}+\ldots
+{x_{m}}^{2}-{x_{m+1}}^{2}-\ldots -x_{2m}^{2}.\notag
\end{equation}
That is, if we write $\varphi = ( \varphi^{1}, \ldots, \varphi^{n+1})$, then\\
$\varphi^{1} = F,\; \varphi^{i} =F \circ G_{i} \;(i = 2,\ldots,n+1)$ for $G_{i}
\in O({\Bbb R}^{2km(n)})$.
\end{theorem}

\begin{proof}
Note that if $G_{1}, \ldots, G_{k}\in O({\Bbb R}^{2m(n)})$ then $G_{1}\oplus
\ldots \oplus G_{k}\in O({\Bbb R}^{2km(n)})$. On the other hand, we have seen
from Theorem \ref{EG} that a \qh \;exists precisely in the case $\varphi:{\Bbb
R}^{2km(n)}\longrightarrow {\Bbb R}^{n+1}$, where we have
\begin{equation}
\varphi \overset{d}{\sim} \lambda_{1}\varphi_{1}\oplus \ldots \oplus
\lambda_{k}\varphi_{k},\notag
\end{equation}
for domain-minimal $\varphi_{i} \in H_{2}^{1}(2m(n),n+1)$, given by irreducible
Clifford systems. Thus it suffices to show the following\\
{\bf Claim:} Any domain-minimal $\varphi_{i} \in H_{2}^{1}(2m(n),n+1)$ arises
from a single quadratic function $F_{0}$.\\
{\bf Proof of Claim:} It follows from Corollary \ref{4.7} and the
Classification Theorem that
\begin{equation}
\varphi_{i} \overset{d}{\sim} \left( X^{t}\left(
                                 \begin{array}{cc}
                                   I_{m}  & 0 \\
                                   0 & -I_{m}
                                  \end{array}
                              \right)
                                   X,\; X^{t}
                              A_{1}X, \ldots,X^{t}
                              A_{n}X\right).\notag
\end{equation}
where the component matrices $A_{i}$, by the Rank Lemma in \cite{OuWoo95}, have
the same rank, the same index and the same spectrum as $\left(
\begin{array}{cc} I_{m}  & 0 \\ 0 & -I_{m} \end{array} \right)$ does. Therefore
from the theory of real quadratic forms we know that there exist $
G_{i}^{\alpha} \in O({\Bbb R}^{2m(n)}) $ such that  $\varphi_{i}^{\alpha} =
\varphi_{i}^{1}\circ G_{i}^{\alpha} = F_{0}\circ G_{i}^{\alpha}\;(\alpha =
2,\ldots,n+1)$, which ends the proof of the claim.
\end{proof}
\begin{example}\label{Ex1}
We can check that ${\varphi}: {\Bbb R}^{8}\longrightarrow {\Bbb R}^{3}$ given
by
\begin{align}
\varphi = &( 2x_{1}^{2}+2x_{2}^{2} + 3x_{3}^{2}+3x_{4}^{2}
-2x_{5}^{2}-2x_{6}^{2} -3x_{7}^{2}-3x_{8}^{2},\notag\\
          &4x_{1}x_{5} + 4x_{2}x_{6} + 6x_{3}x_{8} - 6x_{4}x_{7},\notag\\
         -&4x_{1}x_{6}+4x_{2}x_{5}+6x_{3}x_{7} + 6x_{4}x_{8} )\notag
\end{align}
is a quadratic harmonic morphism. Let \;$ F =  2F_{0} \oplus 3F_{0}$ \;with
\begin{equation}
F_{0}:{\Bbb R}^{4}\longrightarrow {\Bbb R},\;\;F_{0}(x_{1}, \ldots, x_{4}) =
x_{1}^{2} + x_{2}^{2}-x_{3}^{2}-x_{4}^{2}.\notag
\end{equation}
We can further check that $\varphi$ arises from $F$ since
\begin{align}
&\varphi^{1} = F,\;\;\varphi^{2} = F\circ (G_{1} \oplus G_{2}) =F\circ
G^{1},\notag\\
&\varphi^{3} = F\circ (G_{2} \oplus G_{3}) =F\circ G^{2}.\notag
\end{align}
for
\begin{align}
& G_{1} = \left(\begin{array}{cccc}
 \frac{1}{\sqrt{2}} & 0 & \frac{1}{\sqrt{2}} & 0 \\
 0 & \frac{1}{\sqrt{2}} & 0 & \frac{1}{\sqrt{2}} \\
 -\frac{1}{\sqrt{2}} & 0 & \frac{1}{\sqrt{2}} & 0\\
 0 & -\frac{1}{\sqrt{2}} & 0 & \frac{1}{\sqrt{2}}
\end{array}\right),\;
G_{2} = \left(\begin{array}{cccc}
 \frac{1}{\sqrt{2}} & 0 & 0 & -\frac{1}{\sqrt{2}}\\
 0 & \frac{1}{\sqrt{2}} & \frac{1}{\sqrt{2}} & 0 \\
 0 & -\frac{1}{\sqrt{2}}& \frac{1}{\sqrt{2}} & 0\\
 \frac{1}{\sqrt{2}} & 0 & 0 & \frac{1}{\sqrt{2}}
\end{array}\right),\notag\\
& G_{3} = \left(\begin{array}{cccc}
 \frac{1}{\sqrt{2}} & 0 & 0 \frac{1}{\sqrt{2}} \\
 0 & \frac{1}{\sqrt{2}} & -\frac{1}{\sqrt{2}} & 0 \\
 0 & \frac{1}{\sqrt{2}} & \frac{1}{\sqrt{2}} & 0 \\
 -\frac{1}{\sqrt{2}} & 0 & 0 & \frac{1}{\sqrt{2}}
\end{array}\right) \; \in O({\Bbb R}^{4}).\notag
\end{align}
\end{example}

Ferus, Karcher and M\"{u}nzner have noted in \cite{FerKarMun81} that some
Clifford systems can be extended, by adding one member, to be a Clifford system
with one more range-dimension. Our next theorem gives the conditions on the
{\bf range-extendability} of Clifford systems, O-systems and \qh s.

\begin{theorem}\label{extend}
Any domain-minimal \qh s\; (respectively,Clifford systems, or O-systems) which
are not range-maximal can be extended, by adding component functions
(respectively, system members), to be a range-maximal one.
\end{theorem}

\begin{proof}
By Corollary \ref{4.7}, any domain-minimal \qh \; $\varphi:{\Bbb
R}^{2m(n)}\longrightarrow {\Bbb R}^{n+1}$ is $Q$-nonsingular and umbilical, and
hence it arises from a single quadratic function $\lambda F_{0}$ for
$F_{0}:{\Bbb R}^{2m(n)}\longrightarrow {\Bbb R}$,
\begin{equation}
F_{0}(x_{1}, \ldots, x_{2m}) = {x_{1}}^{2}+\ldots
+{x_{m}}^{2}-{x_{m+1}}^{2}-\ldots -x_{2m}^{2}.\notag
\end{equation}
On the other hand, any domain-minimal and range-maximal \qh \; $\varphi:{\Bbb
R}^{2m}\longrightarrow {\Bbb R}^{\sigma(m)+1}$ also arises from the quadratic
function $\lambda F_{0}$. Thus if $n\leq \sigma(m)$ we can add some component
functions of the form $ \lambda F_{0}\circ G^{\alpha}$ untill we get a
range-maximal one. This proves the result for \qh s. The corresponding results
for Clifford systems and O-systems follow from the relationships (Theorems
\ref{P1} and \ref{OC}) between \qh s, Clifford systems and O-systems.
\end{proof}

Note that the \qh \;in Example \ref{Ex1} is not domain-minimal. Though it is
not range-maximal either, it cannot be extended to be a range-maximal one as
one can check easily. On the other hand, the standard multiplication of complex
numbers is neither domain-minimal nor range-maximal. But, as we will see from
the following remark that it can be extended to be a range-maximal one.

\begin{remark}
Theorem \ref{extend} gives a method of constructing \qh s from some given ones.
It is interesting to note that this construction includes the Hopf construction
maps of the standard multiplication $p_{n}:{\Bbb R}^{2n}\longrightarrow {\Bbb
R}^{n},\;(n=1, 2, 4\; or\; 8)$ of real, complex, quaternionic and Cayley
numbers as special cases: Since $ 2 p_{n}:{\Bbb R}^{2n}\longrightarrow {\Bbb
R}^{n}$ is also a \qh \;and for $n= 4$ or $8$ they are domain-minimal but not
range-maximal. Therefore, by Theorem \ref{extend}, they can be extended, by
adding one component function to be range-maximal one as $H(x,y) = (\|x\|^{2}
-\|y\|^{2},\; 2p_{n}(x,y))$, which is exactly the Hopf construction map. The
cases for $n= 1$ or $2$ are easier to check. It can be further checked that in
the above cases there are at most two possible ways of adding one component
function in doing the extension.
\end{remark}

\begin{theorem}
For \;$ n = 1,\;4$ or $8$, any \qh \; $\varphi \in H_{2}(2n,n)$ is
domain-equivalent to a constant multiple of the standard multiplications of the
real algebras of real, quaternionic and Cayley numbers respectively.
\end{theorem}

\begin{proof}
We first note that in all three cases in question, $\varphi$ is domain-minimal
and hence $Q$-nonsingular and umbilical. Therefore , by $(a)$ of Corollary
\ref{4.5}, $ \varphi \overset{d}{\sim}\lambda \varphi_{0}$ for $\varphi_{0} \in
H_{2}^{1}(2n,n)$ given by an irreducible Clifford system $\{P_{i}\} \in
C(2n,n)$ which, in all three cases, belongs to exactly one algebraically
equivalent class (see \cite{FerKarMun81}). On the other hand, it is known (see
Baird \cite{Bai83} Theorem 7.2.7) that the standard multiplications of the real
algebras of real, quaternionic and Cayley numbers are, respectively, in the
class. Thus we obtain the theorem.
\end{proof}

\begin{proposition}\label{H}
For \;$ n = 1,\;2,\;4$ or $8$, any \qh \;$\varphi \in H_{2}(2n,n+1)$ is
bi-equivalent to a constant multiple of the Hopf construction map in the
corresponding cases, and therefore, $\varphi$ \;restricts to \h \; ${\Bbb
S}^{2n-1}\longrightarrow {\Bbb S}^{n}(\lambda)$, where ${\Bbb S}^{n}(\lambda)$
denotes the Euclidean sphere of radius $\lambda$.
\end{proposition}

\begin{proof}
$n=1$ is trivial. For $ n = 2$ the results have been obtained in
\cite{OuWoo95}, where all \qh s $\varphi \in H_{2}(4,3)$ are determined
explicitly and are proved to be domain-equivalent to a constant multiple of the
Hopf construction map. Now we give a proof which treats all four cases. Note
that in all these cases, any \qh \;$\varphi$ is domain-minimal (see Theorem
\ref{EU}) and hence $Q$-nonsingular and umbilical. Therefore , $\varphi
\overset{d}{\sim} \lambda \varphi_{0}$, for $\varphi_{0}$ given by an
irreducible Clifford system $\{P_{i}\} \in C(2n,n+1)$. Now the first half of
the statement in the proposition follows from Proposition \ref{P1} and the fact
(see \cite{FerKarMun81}) that there exists exactly one geometric equivalence
classes in these cases. The second half of the statement is trivial and is
omitted.
\end{proof}

\begin{remark}
{}From Proposition \ref{H} it follows that any \qh \;$ \varphi:{\Bbb
R}^{2n}\longrightarrow {\Bbb R}^{n+1},\;(n= 1,\;2, 4\; or\; 8)$ restricts to a
composition of the classical Hopf fibration followed by a homothety. Thus, in a
sense, Proposition \ref{H} generalizes part of a result of Eells and Yiu (see
\cite{EelYiu94}) which says that if $ \phi: {\Bbb S}^{m}\longrightarrow {\Bbb
S}^{n}$ is the restriction of a homogeneous polynomial harmonic morphism $\Phi:
{\Bbb R}^{m+1}\longrightarrow {\Bbb R}^{n+1}$. Then $\phi$ is isometric to the
classical Hopf fibration.
\end{remark}

\begin{remark}
Concerning \qh s between spheres we know (see \cite{EelYiu94}) that a quadratic
map $\varphi: {\Bbb S}^{m} \longrightarrow {\Bbb S}^{n}$ is a harmonic morphism
if and only if it is isometric to one of the classical Hopf fibrations ${\Bbb
S}^{2n-1} \longrightarrow {\Bbb S}^{n}$ for $n = 1, 2, 4$, or $8$.
\end{remark}

{}From Remark \ref{3.4} we see that any umbilical $\varphi \in
H_{2}^{1}(2m,n+1)$ is associated with an orthogonal multiplication
$\mu_{\varphi} \in F(n,m;m)$ by

\begin{equation}
\varphi \mapsto \{P_{\alpha}\} \mapsto \{\tau_{i}\} \mapsto
\mu_{\varphi}.\label{m5}
\end{equation}
It is interesting to note that the Hopf construction maps of the standard
multiplications of the real algebras of real, complex, quaternionic and Cayley
numbers are the only umbilical \qh s which correspond to \om s that are also
harmonic morphisms.

\begin{proposition}
Let $\varphi \in H_{2}^{1}(2m,n+1)$, and $\mu_{\varphi} \in F(n,m;m)$ be the
corresponding orthogonal multiplication via (\ref{m5}). Then $\mu_{\varphi}$ is
a \h \; if and only if $\varphi$ is bi-equivalent to the Hopf construction maps
of the standard multiplications of real, complex, quaternionic and Cayley
numbers.
\end{proposition}

\begin{proof}
{}From Baird \cite{Bai83} it follows that, an orthogonal multiplication
\\$\mu_{\varphi} \in F(n,m;m)$ is a \h \; if and only if $ n=m
=1,\;2,\;4,\;or\; 8$. This implies that $\varphi \in H_{2}^{1}(2n,n+1)$ for $
n= 1,\;2,\;4,\;or\; 8$ which, together with Proposition \ref{H}, gives the
required results.
\end{proof}

It is easily checked that any \qh \;${\Bbb R}^{2}\longrightarrow {\Bbb R}^{2}$
is domain-equivalent to $\varphi_{0}(z) = z^{2} :{\Bbb R}^{2}\cong {\Bbb
C}\longrightarrow {\Bbb R}^{2}\cong {\Bbb C}$. A proof similar to that of
Proposition \ref{H} can be applied to give the following

\begin{theorem}\label{TL}
Any \qh \;$ \varphi: {\Bbb R}^{m}\longrightarrow {\Bbb R}^{2}$ is
domain-equivalent to
\begin{equation}
\lambda_{1}z_{1}^{2} + \ldots + \lambda_{k}z_{k}^{2}:\;\;{\Bbb
C}^{k}\longrightarrow {\Bbb C},\notag
\end{equation}
where $k=[\frac{m}{2}]$, and $\lambda_{i} \geq 0 \;(i = 1,\ldots,k)$ with at
least one not zero.
\end{theorem}

\begin{ack} This work was done while the author was visiting Department of Pure
Mathematics, University of Leeds, U.K.. The author is grateful to the
department for the hospitality and generosity, and to many colleagues there for
their help  and friendship. Especially, the author is indebted to John C. Wood
for making the author's visit to Leeds possible and for his continued help
through many conversations during the preparation of this work. The author also
wishes to thank J. Eells for some valuable comments and suggestions related to
this work.
\end{ack}

\end{document}